\documentclass[11pt]{article}
\usepackage[russian]{babel}
\usepackage[cp1251]{inputenc}
\usepackage{amsmath,latexsym}
\usepackage[psamsfonts]{amssymb}
\usepackage[dvips]{graphicx}
\usepackage[mathcal]{euscript}
\begin{document}
	MSC2010.35M10
\begin{center}
\textbf{BOUNDARY VALUE PROBLEM FOR HIGH ORDER EQUATION
WITH DISCONTINUOUS COEFFICIENTS
}\\B.Y.Irgashev
 \\Engineering-Construction Institute in Namangan.\\ Namangan,
Uzbekistan
\\
  E-mail : bahromirgasev@gmail.com\end{center}

\textbf{Abstract:} The article considers the Dirichlet problem for a high-order mixed-type equation that splits into factors, each of which is a Lavrentiev-Bitsadze equation with its own excellent coefficient. Sufficient conditions are found for the coefficients under which the problem has a classical solution.

\textbf{Keywords:} Lavrentiev-Bitsadze equation, Fourier coefficients, uniqueness, existence, "small" denominators, series, convergence, Parseval equality.

\textbf{1.Introduction.} Consider a partial differential equation
\[
Lu \equiv \left\{ {\prod\limits_{j = 1}^n {\left( {a_j^2
\frac{{\partial ^2 }}{{\partial x^2 }} + {\mathop{\rm sgn}}
y\frac{{\partial ^2 }}{{\partial y^2 }}} \right)} }
\right\}u\left( {x,y} \right) = 0,\eqno(1)
\]

in a rectangular area  $\Omega  = \left\{ {\left( {x,y} \right):0
< x < 1, - 1 < y < 1} \right\},$ где $0 < a_i  \in R,\,\,a_i  > a_j$ при $i > j.$
 Пусть $\Omega _ +   = \Omega  \cap \left( {y
> 0} \right),\,\Omega _ - = \Omega  \cap \left( {y < 0} \right).$

\textbf{Dirichlet problem.} Find the function $ u (x, y) $ in the domain $ D $
satisfying the conditions:
\[
u \in C^{2n - 1} \left( {\overline \Omega  } \right) \cap C^{2n}
\left( {\Omega _ +   \cup \Omega _ -  } \right),\eqno(2)
\]
\[
Lu\left( {x,y} \right) \equiv 0,\,\,\,\left( {x,y} \right) \in
\Omega _ +   \cup \Omega _ -  ,\eqno(3)
\]
\[
D_x^{2s} u\left( {0,y} \right) = D_x^{2s} u\left( {1,y} \right) =
0,\,\,\, - 1 \le y \le 1,\,\eqno(4)
\]
\[D_y^{2s}u\left( {x,1} \right) = {\varphi _s}\left( x \right),\,\,0 \le x \le 1,\eqno(5)
\]
\[
D_y^{2s} u\left( {x, - 1} \right) = \psi _s \left( x \right),\,\,0
\le x \le 1,\eqno(6)
\]

where $ s = 0,n - 1,\,\,\,\,\,\varphi_s\left(x\right),\psi _s
\left(x\right)$ -sufficiently smooth and terms of agreement works for them.
Equation in the form
\[
Lu \equiv \left\{ {\prod\limits_{j = 1}^n {\left( {\frac{{\partial ^2 }}{{\partial t^2 }} - a_j^2 \Delta  + c_j } \right)} } \right\}u\left( {t,x} \right) = f\left( {t,x} \right),
\]
where
\[
\Delta  = \frac{{\partial ^2 }}{{\partial x_1^2 }} + \frac{{\partial ^2 }}{{\partial x_2^2 }} + ... + \frac{{\partial ^2 }}{{\partial x_m^2 }},\,\,a_j ,c_j  \in R,\,\,j = 1,2,...,n,
\]

are found in the classical theory of elasticity [1], as well as in quantum theory field [2].
The importance of studying equations of type (1) was noted by A.V. Bitsadze [3].
Equation (1), in which all coefficients are equal, were studied in the works of M.M.Smirnov [4], [5], V.I. Zhegalov [6], and K. B. Sabitova [7], [8]. We also note [9], in which operator equations were considered. High-order equations with smooth coefficients, such as equation (1), were studied in [10] - [13] and others.
For $n =1$, equation (1) is the well-known Lavrentiev–Bitsadze equation for which A.V.Bitsadze showed the incorrectness of the Dirichlet problem [14]. Recently, many specialists of the Lavrentiev-Bitsadze type equation, using the spectral method, have studied various boundary value problems [15], [16] and others. In these works ([15], [16] and others) sufficient conditions were obtained for the coefficient ${a_1}$ for which the boundary value problems under consideration have a classical solution.

In this paper, a criterion for the uniqueness of a solution to the Dirichlet problem is established. The solution is constructed as a sum over the eigenfunctions of the one-dimensional problem. When substantiating the convergence of a series, the problem of “small” denominators arises. Sufficient conditions are obtained for the coefficients ${a_j}\,\left({j = 1, ...,n}\right)\,$, for which exists a classical solution of the Dirichlet problem.

\textbf{2.The uniqueness of the solution.} Let $u(x,y)$ be a solution of
equations (1) with conditions (2) - (6). Consider its Fourier coefficients
\[
u_k \left( y \right) = \sqrt 2 \int\limits_0^1 {u\left( {x,y}
\right)\sin \pi kxdx} ,\,\,k = 1,2,...\,.\eqno(7)
\]

Based on (7), we introduce the functions
\[
u_{k,\varepsilon } \left( y \right) = \sqrt 2
\int\limits_\varepsilon ^{1 - \varepsilon } {u\left( {x,y}
\right)\sin \pi kxdx},\eqno(8)
\]
where $\varepsilon>0$ is a fairly small number. Differentiating
equality (8) with respect to $y$ under the sign of the integral $2n$ times, for $y>0$ and
$y<0$, given equation (1) and passing to the limit as
$\varepsilon\rightarrow+0$, taking into account the boundary conditions (4) - (6),
the continuity of the function $u(x,y)$ and its derivatives on the line $y=0$,
we get the following task:
\[\left\{ \begin{array}{l}
\prod\limits_{j = 1}^n {\left( { - {{\left( {{a_j}k\pi } \right)}^2} + {\mathop{\rm sgn}} y\frac{{{\partial ^2}}}{{\partial {y^2}}}} \right)u\left( {x,y} \right)}  = 0,\\
u_k^{2s}\left( 1 \right) = {\varphi _{sk}},\\
u_k^{2s}\left( { - 1} \right) = {\psi _{sk}},\\
u_k^{\left( t \right)}\left( { + 0} \right) = u_k^{\left( t \right)}\left( { - 0} \right),\,s = \overline {0,n - 1} ,\,\,\,t = \overline {0,2n - 1} .
\end{array} \right.\eqno(9)
\]
where
\[
\varphi _{sk}  = \sqrt 2 \int\limits_0^1 {\varphi _s \left( x
\right)\sin \pi kxdx,\,\,} \psi _{sk}  = \sqrt 2 \int\limits_0^1
{\psi _s \left( x \right)\sin \pi kxdx\,\,} .
\]

Then, for $y>0$, the general solution of equation (9) has the form
\[
u_k \left( y \right) = \sum\limits_{j = 1}^n {\left( {c_{2j - 1}
e^{a_j \pi ky}  + c_{2j} e^{ - a_j \pi ky} } \right)} ,
\]
\[u_k^{\left( {2m} \right)}\left( y \right) = {\left( {\pi k} \right)^{2m}}\sum\limits_{j = 1}^n {\left( {{c_{2j - 1}}a_j^{2m}{e^{{a_j}\pi ky}} + {c_{2s}}a_j^{2m}{e^{ - {a_j}\pi ky}}} \right)} ,\]
where
\[
m = 0,...,2n - 1.
\]
For $y<0$
\[
u_k \left( y \right) = \sum\limits_{s = 1}^n {\left( {d_{2s - 1}
\cos \pi ka_s y + d_{2s} \sin \pi ka_s y} \right)} ,
\]
\[u_k^{\left( {2m} \right)} = {\left( {\pi k} \right)^{2m}}\sum\limits_{s = 1}^n {\left( {{d_{2s - 1}}{{\left( { - 1} \right)}^m}a_s^{2m}\cos \left( {\pi k{a_s}y} \right) + {d_{2s}}{{\left( { - 1} \right)}^m}a_s^{2m}\sin \left( {\pi k{a_s}y} \right)} \right)} .\]

Satisfying the boundary conditions, we obtain the system
\[\left\{ \begin{array}{l}
\sum\limits_{s = 1}^n {\left( {{c_{2s - 1}}a_s^{2m}{e^{\pi k{a_s}}} + {c_{2s}}a_s^{2m}{e^{ - \pi k{a_s}}}} \right)}  = \frac{{{\varphi _{mk}}}}{{{{\left( {\pi k} \right)}^{2m}}}},\,\,\\
\sum\limits_{s = 1}^n {\left( {{d_{2s - 1}}a_s^{2m}\cos \pi k\left( { - {a_s}} \right) + {d_{2s}}a_s^{2m}\sin \pi k\left( { - {a_s}} \right)} \right)}  = {\left( { - 1} \right)^m}\frac{{{\psi _{mk}}}}{{{{\left( {\pi k} \right)}^{2m}}}},\\
m = 0,...,n - 1,\\
\sum\limits_{s = 1}^n {\left( {{c_{2s - 1}}a_s^t + {c_{2s}}{{\left( { - {a_s}} \right)}^t}} \right)}  - \sum\limits_{s = 1}^n {\left( {{d_{2s - 1}}a_s^t\cos \frac{{\pi t}}{2} + {d_{2s}}a_s^t\sin \frac{{\pi t}}{2}} \right) = 0} ,\\
t = 0,...,2n - 1.
\end{array} \right.\eqno(10)
\]
The structure of the main determinant of system (10) has the form:
\[
\Delta  = \Delta _1 \left( {\Delta _2  + \Delta _k } \right),
\]
where
\[{\Delta _1} = {e^{\pi k\left( {{a_1} + ... + {a_n}} \right)}}\left| {\begin{array}{*{20}{c}}
1&1&.&.&1\\
{a_1^2}&{a_2^2}&.&.&{a_n^2}\\
{a_1^4}&{a_2^4}&.&.&.\\
.&.&.&.&.\\
{a_1^{2(n - 1)}}&{a_2^{2(n - 1)}}&.&.&{a_n^{2(n - 1)}}
\end{array}} \right| = {e^{\pi k\left( {{a_1} + ... + {a_n}} \right)}}\prod\limits_{j > s} {\left( {a_j^2 - a_s^2} \right)} ,\]
\[
\Delta _2  = \left( { - \frac{i}{2}} \right)^n \left|
{\begin{array}{*{20}c}
   \begin{array}{l}
 0... \\
 . \\
 \end{array} & \begin{array}{l}
 0 \\
 . \\
 \end{array} & \begin{array}{l}
 e^{ - i\pi ka_1 }  \\
 . \\
 \end{array} & \begin{array}{l}
 e^{i\pi ka_1 } ... \\
 . \\
 \end{array} & \begin{array}{l}
 e^{ - i\pi ka_n }  \\
 . \\
 \end{array} & \begin{array}{l}
 e^{i\pi ka_n }  \\
 . \\
 \end{array}  \\
   {0...} & 0 & {a_1^{n - 1} e^{^{ - i\pi ka_1 } } } & {a_1^{n - 1} e^{i\pi ka_1 } ...} & {a_n^{n - 1} e^{ - i\pi ka_n } } & {a_n^{n - 1} e^{i\pi ka_n } }  \\
   {1...} & 1 & 1 & 1 & 1 & 1  \\
   { - a_1 ...} & { - a_n } & {ia_1 } & { - ia_1 ...} & {ia_n } & { - ia_n }  \\
   . & . & . & . & . & .  \\
   {\left( { - a_1 } \right)^{2n - 1} ...} & {\left( { - a_n } \right)^{2n - 1} } & {\left( {ia_1 } \right)^{2n - 1} } & {\left( { - ia_1 } \right)^{2n - 1} ...} & {\left( {ia_n } \right)^{2n - 1} } & {\left( { - ia_n } \right)^{2n - 1} }  \\
\end{array}} \right|,
\]
\[
\mathop {\lim }\limits_{k \to  + \infty } \Delta _k  = 0.
\]

\textbf{Theorem 1.} Е\emph{If there is a solution to the Dirimhlet problem, then it is unique only if the condition $\Delta\ne0$ is satisfied, for all $k$}

	\textbf{Proof.} Suppose that $\Delta\ne0$ and boundary conditions (4) - (6) are homogeneous, while system (10) has only the trivial solution ${u_k}\left(y\right)\\equiv 0,$for all numbers $k$. So $u\left({x,y}\right)=0$ almost everywhere, but since $u\left({x, y}\right)$ is continuous in $\overline\Omega $, then  $u\left( {x,y} \right) \equiv 0$ in$\overline\Omega$.
   Now, for some values ??of $k$, let the main determinant of $\Delta=0$.Then the homogeneous Dirichlet problem will have a nonzero solution, which will violate the uniqueness of the solution.
\textbf{Theorem is proved.}

\textbf{3.Existence of solution.}
Having made some transformations, we have

\[
\Delta _2  = C_1 i^n \left| {\begin{array}{*{20}c}
   \begin{array}{l}
 0... \\
 0... \\
 . \\
 \end{array} & \begin{array}{l}
 0 \\
 0 \\
 . \\
 \end{array} & \begin{array}{l}
 e^{ - i\pi ka_1 }  \\
 0 \\
 . \\
 \end{array} & \begin{array}{l}
 e^{i\pi ka_1 } ... \\
 0... \\
 . \\
 \end{array} & \begin{array}{l}
 0 \\
 0 \\
 . \\
 \end{array} & \begin{array}{l}
 0 \\
 0 \\
 . \\
 \end{array}  \\
   {0...} & 0 & 0 & {0...} & {e^{ - i\pi ka_n } } & {e^{i\pi ka_n } }  \\
   {1...} & 1 & 1 & 1 & 1 & 1  \\
   {a_1 ...} & {a_n } & { - ia_1 } & {ia_1 ...} & { - ia_n } & {ia_n }  \\
   . & . & . & . & . & .  \\
   {a_1 ^{2n - 1} ...} & {a_n ^{2n - 1} } & {\left( { - ia_1 } \right)^{2n - 1} } & {ia_1 ^{2n - 1} ...} & {\left( { - ia_n } \right)^{2n - 1} } & {ia_n ^{2n - 1} }  \\
\end{array}} \right|,
\]

where $0\ne C_1$ is a certain number which is independent from $k.$ \\
Now we state the following theorem, which gives an estimate for the function ${u_k}\left(y\right)$ and its derivatives.

\textbf{Theorem 2.}\emph {To solve ${u_k}\left(y\right)$  problem (9), for $y\ne0$ and sufficiently large values of $k$, estimates
}
\[u_k^{\left( j \right)}\left( y \right) \le M{k^j}\frac{{\sum\limits_{s = 0}^{n - 1} {\left( {\left| {{\varphi _{sk}}} \right| + \left| {{\psi _{sk}}} \right|} \right)} }}{{\left| {{\Delta _2} + {\Delta _k}} \right|}},\eqno(11)
\]
where
\[0 < M = const,\,j = 0,1,...,2n - 1.\]

\textbf{Proof.} Assume,  $y<0$  ( case $y>0$
considered similarly)and all constants will be denoted by a single letter $M$. From the representation of the solution to problem (9) we obtain
\[\left| {u_k^{\left( j \right)}(y)} \right| \le M{k^j}\left| {{u_k}\left( y \right)} \right|,\]
hence it is sufficient to prove the estimate in cases $j=0$. We have
\[\left| {{u_k}\left( y \right)} \right| \le \sum\limits_{s = 1}^n {\left( {\left| {{d_{2s - 1}}} \right| + \left| {{d_{2s}}} \right|} \right)} ,\]
$$\left| {{d_{2s - 1}}} \right| = \left| {\frac{{\Delta \left( {{d_{2s - 1}}} \right)}}{\Delta }} \right| \le M\frac{{\sum\limits_{s = 1}^n {\left\{ {\left| {{\varphi _{ks}}} \right| + \left| {{\psi _{ks}}} \right|} \right\}} \, \cdot O\left( {\left| {{\Delta _1}} \right|} \right)}}{{\left| {{\Delta _1}\left( {{\Delta _2} + {\Delta _k}} \right)} \right|}} \le M\frac{{\sum\limits_{s = 1}^n {\left\{ {\left| {{\varphi _{ks}}} \right| + \left| {{\psi _{ks}}} \right|} \right\}} \,}}{{\left| {{\Delta _2} + {\Delta _k}} \right|}},$$
similarly
$$\left| {{d_{2s}}} \right| = \left| {\frac{{\Delta \left( {{d_{2s}}} \right)}}{\Delta }} \right| \le M\frac{{\sum\limits_{s = 1}^n {\left\{ {\left| {{\varphi _{ks}}} \right| + \left| {{\psi _{ks}}} \right|} \right\}} \,O\left( {\left| {{\Delta _1}} \right|} \right)}}{{\left| {{\Delta _1}\left( {{\Delta _2} + {\Delta _k}} \right)} \right|}} \le M\frac{{\sum\limits_{s = 1}^n {\left\{ {\left| {{\varphi _{ks}}} \right| + \left| {{\psi _{ks}}} \right|} \right\}} \,}}{{\left| {{\Delta _2} + {\Delta _k}} \right|}},$$
 here $\Delta \left( {{d_{2s - 1}}} \right),\Delta \left( {{d_{2s}}} \right)$
determinants obtained from the determinant $\Delta$ by replacing the corresponding columns with the right-hand side of system (10).

 \textbf{Theorem is proved.}

 As can be seen from (11), the “small denominator” problem may arise. Let us find out the conditions for the separation from zero of the denominator in expression (11). To do this, it is enough to obtain the conditions for the ${\Delta_2}.$ expression to be separated from zero. Consider $\Delta _2$, for some particular values ??of $a_j$.
\\
\textbf{1.  Let $a_j  \in N,$  for $j=1,2,...,n.$} Then
\[
\Delta _2  = C_1 i^n \left| {\begin{array}{*{20}c}
   0 & {...} & 0 & {\left( { - 1} \right)^{ka_1 } } & {\left( { - 1} \right)^{ka_1 } } & {...} & 0 & 0  \\
   0 & {...} & 0 & 0 & 0 & {...} & 0 & 0  \\
   . & . & . & . & . & . & . & .  \\
   0 & {...} & . & . & . & {...} & {\left( { - 1} \right)^{ka_n } } & {\left( { - 1} \right)^{ka_n } }  \\
   1 & {...} & 1 & 1 & 1 & {...} & 1 & 1  \\
   {a_1 } & {...} & {a_n } & { - ia_1 } & {ia_1 } & {...} & { - ia_n } & {ia_n }  \\
   . & . & . & . & . & . & . & .  \\
   {a_1^{2n - 1} } & {...} & {a_n^{2n - 1} } & {\left( { - ia_1 } \right)^{2n - 1} } & {\left( {ia_1 } \right)^{2n - 1} } & . & {\left( { - ia_n } \right)^{2n - 1} } & {\left( {ia_n } \right)^{2n - 1} }  \\
\end{array}} \right| =
\]
\[
 = C_1 i^n \left| {\begin{array}{*{20}c}
   0 & {...} & 0 & {\left( { - 1} \right)^{ka_1 } } & 0 & {...} & 0 & 0  \\
   0 & {...} & 0 & 0 & 0 & {...} & 0 & 0  \\
   . & . & . & . & . & . & . & .  \\
   0 & {...} & . & . & . & {...} & {\left( { - 1} \right)^{ka_n } } & {\left( { - 1} \right)^{ka_n } }  \\
   1 & {...} & 1 & 1 & 0 & {...} & 1 & 0  \\
   {a_1 } & {...} & {a_n } & { - ia_1 } & {2ia_1 } & {...} & { - ia_n } & {2ia_n }  \\
   . & . & . & . & . & . & . & .  \\
   {a_1^{2n - 1} } & {...} & {a_n^{2n - 1} } & {\left( { - ia_1 } \right)^{2n - 1} } & {2\left( {ia_1 } \right)^{2n - 1} } & . & {\left( { - ia_n } \right)^{2n - 1} } & {2\left( {ia_n } \right)^{2n - 1} }  \\
\end{array}} \right| =
\]
\[
 = C_2 \left| {\begin{array}{*{20}c}
   1 & {...} & 1 & 0 & {...} & 0  \\
   \begin{array}{l}
 a_1^2  \\
 . \\
 \end{array} & \begin{array}{l}
 ... \\
 . \\
 \end{array} & \begin{array}{l}
 a_n^2  \\
 . \\
 \end{array} & \begin{array}{l}
 0 \\
 . \\
 \end{array} & \begin{array}{l}
 ... \\
 . \\
 \end{array} & \begin{array}{l}
 0 \\
 . \\
 \end{array}  \\
   {a_1^{2n - 2} } & {...} & {a_n^{2n - 2} } & 0 & {...} & 0  \\
   {a_1 } & {...} & {a_n } & {2ia_1 } & {...} & {2ia_n }  \\
   . & . & . & . & . & .  \\
   {a_1^{2n - 1} } & {...} & {a_n^{2n - 1} } & {2\left( {ia_1 } \right)^{2n - 1} } & {...} & {2\left( {ia_n } \right)^{2n - 1} }  \\
\end{array}} \right| =
\]
\[
 = C_2 \left| {\begin{array}{*{20}c}
   1 & {...} & 1  \\
   {a_1^2 } & {...} & {a_n^2 }  \\
   . & . & .  \\
   {a_1^{2n - 2} } & {...} & {a_n^{2n - 1} }  \\
\end{array}} \right|\left| {\begin{array}{*{20}c}
   {2ia_1 } & {...} & {2ia_n }  \\
   {2\left( {ia_1 } \right)^3 } & {...} & {2\left( {ia_n } \right)^3 }  \\
   . & . & .  \\
   {2\left( {ia_1 } \right)^{2n - 1} } & {...} & {2\left( {ia_n } \right)^{2n - 1} }  \\
\end{array}} \right| \ne 0,
\]
here $C_{2}$ - non-zero constant, independent from $k$.

\textbf{Conclusion:}  $\Delta_{2}$,up to a sign, does not depend on $k$ and is nonzero, for ${a_j} \in N,\,j = 1,2,...,n.$

\textbf{2. Let now $ a_j  \in Q,\,\,\forall j$,} then we have
\[
\begin{array}{l}
 \Delta _2  = C_1 i^n e^{ - i\pi k\left( {a_1  + a_2  + ... + a_n } \right)}  \\
 \left| {\begin{array}{*{20}c}
   1 & 1 & . & 1 & {1 - e^{2i\pi ka_1 } } & . & . & {1 - e^{2i\pi ka_n } }  \\
   {a_1 } & {a_2 } & . & {a_n } & {ia_1  - \left( { - i} \right)a_1 e^{2i\pi ka_1 } } & . & . & {ia_n  - \left( { - ia_n } \right)e^{2i\pi ka_n } }  \\
   . & . & . & . & . & . & . & .  \\
   {a_1^{2n - 1} } & {a_2^{2n - 1} } & . & {a_n^{2n - 1} } & {\left( {ia_1 } \right)^{2n - 1}  - \left( { - ia_1 } \right)^{2n - 1} e^{2i\pi ka_1 } } & . & . & {\left( {ia_n } \right)^{2n - 1}  - \left( { - ia_n } \right)^{2n - 1} e^{2i\pi ka_n } }  \\
\end{array}} \right|, \\
 \end{array}
\]
using the following qualifier property
\[
\left| {\begin{array}{*{20}c}
   {a + b} & {c + d}  \\
   {k + l} & {n + m}  \\
\end{array}} \right| = \left| {\begin{array}{*{20}c}
   a & c  \\
   k & n  \\
\end{array}} \right| + \left| {\begin{array}{*{20}c}
   a & d  \\
   k & m  \\
\end{array}} \right| + \left| {\begin{array}{*{20}c}
   b & c  \\
   l & n  \\
\end{array}} \right| + \left| {\begin{array}{*{20}c}
   b & d  \\
   l & m  \\
\end{array}} \right|,
\]
then we have
\[
\Delta _2  = C_1 i^n e^{ - i\pi k\left( {a_1  + a_2  + ... + a_n } \right)} \Delta _3 ,
\]
where
\[
\begin{array}{l}
 \Delta _3  = \left( {a_1 ,...,a_n ,ia_1 ,...,ia_n } \right) - \left( {a_1 ,...,a_n , - ia_1 ,ia_2 ,...,ia_n } \right)e^{2i\pi ka_1 }  -  \\
  - \left( {a_1 ,...,a_n ,ia_1 , - ia_2 ,ia_3 ,...,ia_n } \right)e^{2i\pi ka_2 }  - ... \\
  - \left( {a_1 ,...,a_n ,ia_1 ,..., - ia_n } \right)e^{2i\pi ka_n }  +  \\
  + \left( {a_1 ,...,a_n , - ia_1 , - ia_2 ,...,ia_n } \right)e^{2i\pi k\left( {a_1  + a_2 } \right)}  + ... +  \\
  + \left( {a_1 ,...,a_n , - ia_1 ,ia_2 ,...,ia_{n - 1} , - ia_n } \right)e^{2i\pi k\left( {a_1  + a_n } \right)}  +  \\
  + \left( {a_1 ,...,a_n ,ia_1 , - ia_2 , - ia_3 ,ia_4 ...,ia_n } \right)e^{2i\pi k\left( {a_2  + a_3 } \right)}  + ... +  \\
  + \left( {a_1 ,...,a_n ,ia_1 , - ia_2 ,ia_3 ,...,ia_{n - 1} , - ia_n } \right)e^{2i\pi k\left( {a_2  + a_n } \right)}  \\
 ... + \left( { - 1} \right)^n \left( {a_1 ,...,a_n , - ia_1 ,..., - ia_n } \right)e^{2i\pi k\left( {a_1  + a_2  + ... + a_n } \right)} , \\
 \end{array}\eqno(12)
\]
here
\[
\left( {b_1 ,b_2 ,...,b_n } \right) = \prod\limits_{j > i} {\left( {b_j  - b_i } \right)} .
\]
Let
\[
a_j  = \frac{{s_j }}{{t_j }},\,\,\,a_i  + a_j  = \frac{{s_{ij} }}{{t_{ij} }},...,\,\,\,a_1  + a_2  + ... + a_n  = \frac{{s_{12...n} }}{{t_{12...n} }},
\]
where
\[
\begin{array}{l}
 s_1 ,s_2 ,...,s_{12...n} ,t_1 ,t_2 ,...t_{12...n}  \in N,\,\,\left( {s_j ,t_j } \right) = 1,\,j = 1,...,n; \\
 \left( {s_{ij} ,t_{ij} } \right) = 1,\,i = 1,...,n,\,j = 2,...,n,\,i \ne j;...,\,\left( {s_{12...n} ,t_{12...n} } \right) = 1, \\
 \end{array}
\]
then
\[
\begin{array}{l}
 ka_j  = k_j  + \frac{{m_j }}{{t_j }},\,\,k\left( {a_i  + a_j } \right) = k_{ij}  + \frac{{m_{ij} }}{{t_{ij} }},k\left( {a_i  + a_j  + a_l } \right) = k_{ijl}  + \frac{{m_{ijl} }}{{t_{ijl} }},...,\, \\
 k\left( {a_1  + a_2  + ... + a_n } \right) = k_{12..n}  + \frac{{m_{12...n} }}{{t_{12...n} }}, \\
 \end{array}
\]
where
\[
\begin{array}{l}
 k_1 ,k_2 ,...k_{12...n} ,m_1 ,m_2 ,...,m_{12...n}  \in N;\,\, \\
 0 \le k_i  < t_i ,\,i = 1,...,n;\, \\
 0 \le m_{ij}  < t_{ij} ,\,i = 1,...,n,\,j = 2,...,n,\,i \ne j; \\
 \end{array}
\]
\[
\begin{array}{l}
 0 \le m_{ijl}  < t_{ijl} ,\,i = 1,...,n,\,j = 2,...,n,l = 3,...,n,\,i \ne j \ne l; \\
 0 \le m_{12...n}  < t_{12...n}  = 1. \\
 \end{array}
\]
If all ${k_i},{k_{ij}},{k_{ijl}},...,{k_{12...n}}$  are equal to zero, then we come to the first case, therefore, we assume that there are nonzero integer parts of the fractions $k{a_1},k{a_2}, ...,k\left({{a_1} + ... +{ a_n}}\right)$. Now we transform expression (12) to the form
\[
\begin{array}{l}
 \Delta _3 \left( {m_1 ,m_2 ,...,m_{12...n} } \right) = \left( {a_1 ,...,a_n ,ia_1 ,...,ia_n } \right) - \left( {a_1 ,...,a_n , - ia_1 ,ia_2 ,...,ia_n } \right)e^{2i\pi \frac{{m_1 }}{{t_1 }}}  -  \\
  - \left( {a_1 ,...,a_n ,ia_1 , - ia_2 ,ia_3 ,...,ia_n } \right)e^{2i\pi \frac{{m_2 }}{{t_2 }}}  - ... \\
 \end{array}\eqno(13)
\]
\[
 - \left( {a_1 ,...,a_n ,ia_1 ,..., - ia_n } \right)e^{2i\pi \frac{{m_n }}{{t_n }}}  + \left( {a_1 ,...,a_n , - ia_1 , - ia_2 ,...,ia_n } \right)e^{2i\pi \frac{{m_{12} }}{{t_{12} }}}  + ...
\]
\[
 + \left( {a_1 ,...,a_n , - ia_1 ,ia_2 ,...,ia_{n - 1} , - ia_n } \right)e^{2i\pi \frac{{m_{1n} }}{{t_{1n} }}}  + \left( {a_1 ,...,a_n ,ia_1 , - ia_2 , - ia_3 ,ia_4 ...,ia_n } \right)e^{2i\pi \frac{{m_{23} }}{{t_{23} }}}  + ... +
\]
\[
 + \left( {a_1 ,...,a_n ,ia_1 , - ia_2 ,ia_3 ,...,ia_{n - 1} , - ia_n } \right)e^{2i\pi \frac{{m_{2n} }}{{t_{2n} }}} ... + \left( { - 1} \right)^n \left( {a_1 ,...,a_n , - ia_1 ,..., - ia_n } \right)e^{2i\pi \frac{{m_{12...n} }}{{t_{12...n} }}},
\]
The number of expressions of the form (13) is finite, therefore, among them one can choose the minimum value. But the minimum value should be nonzero, therefore, we should exclude the roots of expressions (13) from consideration. We introduce the notation
\[
z = e^{\frac{{2i\pi }}{M}} ,
\]
where
\[
M = HOK\left( {t_i ,t_{ij} ,t_{ijl} ,...,t_{12...n} } \right),
\]
then from (13) we obtain
\[
\Delta _3 \left( {m_1 ,m_2 ,...,m_{12...n} } \right) = \left( {a_1 ,...,a_n ,ia_1 ,...,ia_n } \right) - \left( {a_1 ,...,a_n , - ia_1 ,ia_2 ,...,ia_n } \right)z^{m_1 T_1 }  -
\]
\[
 - \left( {a_1 ,...,a_n ,ia_1 , - ia_2 ,ia_3 ,...,ia_n } \right)z^{m_2 T_2 }  - ...\eqno(14)
\]
\[
 - \left( {a_1 ,...,a_n ,ia_1 ,..., - ia_n } \right)z^{m_n T_n }  + \left( {a_1 ,...,a_n , - ia_1 , - ia_2 ,...,ia_n } \right)z^{m_{12} T_{12} }  + ...
\]
\[
 + \left( {a_1 ,...,a_n , - ia_1 ,ia_2 ,...,ia_{n - 1} , - ia_n } \right)z^{m_{1n} T_{1n} }  + \left( {a_1 ,...,a_n ,ia_1 , - ia_2 , - ia_3 ,ia_4 ...,ia_n } \right)z^{m_{23} T_{23} }  + ...
\]
\[
 + \left( {a_1 ,...,a_n ,ia_1 , - ia_2 ,ia_3 ,...,ia_{n - 1} , - ia_n } \right)z^{m_{2n} T_{2n} } ... + \left( { - 1} \right)^n \left( {a_1 ,...,a_n , - ia_1 ,..., - ia_n } \right)z^{m_{12...n} T_{12...n} } ,
\]
here
\[
\begin{array}{l}
 T_i  = \frac{M}{{t_i }},\,T_{ij}  = \frac{M}{{t_{ij} }},T_{ijl}  = \frac{M}{{t_{ijl} }},...,T_{12...n}  = \frac{M}{{t_{12...n} }}, \\
 \,i = 1,...,n;\,j = 2,...,n;\,l = 3,...,n;i \ne j \ne l. \\
 \end{array}
\]
the roots of expression (14) have form:
\[{z_j} = {A_j}\left( {{m_1},{m_2},...,{m_{12...n}}} \right) + i{B_j}\left( {{m_1},{m_2},...,{m_{12...n}}} \right),\,j = 1,...,V\left( {{m_1},{m_2},...,{m_{12...n}}} \right),\]
where
\[V\left( {{m_1},{m_2},...,{m_{12...n}}} \right) = \max \left( {{m_1}{T_1},{m_2}{T_2},...,{m_{12...n}}{T_{12...n}}} \right),\,\,{A_j},{B_j} \in R.\]
Hence
\[{e^{\frac{{2i\pi }}{M}}} = \sqrt {A_j^2 + B_j^2} {e^{i\arg {z_j}}},\eqno(15)\]
for the solvability of equation (15), we require
\[\sqrt {A_j^2 + B_j^2}  = 1,\,\,{z_j} = \pi \frac{{{b_j}\left( {{m_1},{m_2},...,{m_{12...n}}} \right)}}{{{c_j}\left( {{m_1},{m_2},...,{m_{12...n}}} \right)}},\,{b_j},{c_j} \in N,\]
then
\[\frac{2}{M} = \frac{{{b_j}}}{{{c_j}}} \Rightarrow \frac{1}{M} = \frac{{{b_j}}}{{2{c_j}}},\eqno(16)\]
if
\[HOD\left( {M,2{c_j}} \right) = 1,\,\,j = 1,...,V,\eqno(17)\]
then equality (16) will not be satisfied and ${\Delta _3}\left( {{m_1},{m_2},...,{m_{12...n}}} \right) \ne 0,$ that is
\[\left| {{\Delta _3}} \right| \ge \mathop {\min }\limits_{{m_1},{m_2},...,{m_{12...n}}} \left| {{\Delta _3}\left( {{m_1},{m_2},...,{m_{12...n}}} \right)} \right| > 0.\]

\textbf{3.Let ${a_j}$-irrational number}, при $j = 1,...,n.$. We obtain the explicit form of the expression ${\Delta _2}.$ The calculations will be performed this time somewhat differently than in the first two cases. The determinant $\ Delta_{2}$ accurate to the sign has the form
\[
\Delta _2  = C_1 \sum\limits_{t_j  =  \pm 1} {\left( { - 1}
\right)^l A\left( {t_1 a_1 ,t_2 a_2 ,...,t_n a_n } \right)},\eqno(18)
\]
where
\[
\begin{array}{l}
 A\left( {t_1 a_1 ,t_2 a_2 ,...,t_n a_n } \right) = i^n \left( {e^{ - i\pi k\left( {t_1 a_1  + t_2 a_2  + ... + t_n a_n } \right)} \left( {a_1 ,a_{2,} ...,a_n ,it_1 a_1 ,it_2 a_2 ,...,it_n a_n } \right) + } \right. \\
  + \left. {\left( { - 1} \right)^n e^{i\pi k\left( {t_1 a_1  + t_2 a_2  + ... + t_n a_n } \right)} \left( {a_1 ,a_{2,} ...,a_n , - it_1 a_1 , - it_2 a_2 ,..., - it_n a_n } \right)} \right), \\
 \end{array}
\]
\[
t_1  = 1,\,\,\,t_j  =  \pm 1,\,j = 2,...,n,
\]
 $l$- quantity of positive $t_{j}$, the sum of $\Sigma$ considered by every possible
 values of $t_{j}$,
\[
\left( {a_1 ,a_{2,} ...,a_n ,it_1 a_1 ,it_2 a_2 ,...,it_n a_n }
\right) = \prod\limits_{j = 2}^n {\prod\limits_{s = 1}^{j - 1}
{\left( {a_j  - a_s } \right)} } \prod\limits_{j = 2}^n
{\prod\limits_{s = 1}^{j - 1} {\left( {it_j a_j  - it_s a_s }
\right)} }  \cdot
\]
\[
\prod\limits_{j > s} {\left( {it_j a_j  - a_s } \right)\left(
{it_s a_s  - a_j } \right)} \prod\limits_{j = 1}^n {\left( {it_j
a_j  - a_j } \right)}  = i^{\frac{{n\left( {n - 1} \right)}}{2}}
\left( { - 1} \right)^n \prod\limits_{j = 1}^n {a_j }
\prod\limits_{j = 2}^n {\prod\limits_{s = 1}^{j - 1} {\left( {a_j
- a_s } \right)} }  \cdot
\]
\[
\prod\limits_{j = 2}^n {\prod\limits_{s = 1}^{j - 1} {\left( {t_j
a_j  - t_s a_s } \right)} } \prod\limits_{j = 1}^n {\left( {1 -
it_j } \right)} \prod\limits_{j > s} {\left\{ {a_j a_s \left( {1 -
t_j t_s } \right) - i\left( {t_j a_j^2  + t_s a_s^2 } \right)}
\right\}}  =
\]
\[
\begin{array}{l}
  = i^{\frac{{n\left( {n - 1} \right)}}{2}} \left( { - 1} \right)^n \sqrt 2 ^n e^{ - \frac{{i\pi \left( {t_1  + t_2  + ... + t_n } \right)}}{4}} \prod\limits_{j = 1}^n {a_j } \prod\limits_{j = 2}^n {\prod\limits_{s = 1}^{j - 1} {\left( {a_j  - a_s } \right)} }  \\
 \prod\limits_{j = 2}^n {\prod\limits_{s = 1}^{j - 1} {\left( {t_j a_j  - t_s a_s } \right)} } \prod\limits_{j > s} {\left( {\left( {a_j^2  + a_s^2 } \right)e^{ - i\alpha _{js} } } \right)} , \\
 \end{array}
\]
where
\[
\sin \alpha _{js}  = \frac{{t_j a_j^2  + t_s a_s^2 }}{{a_j^2  +
a_s^2 }},
\]

here we took advantage of the fact that
\[
\begin{array}{l}
 \left| {a_j a_s \left( {1 - t_j t_s } \right) - i\left( {t_j a_j^2  + t_s a_s^2 } \right)} \right| = \sqrt {\left( {a_j a_s \left( {1 - t_j t_s } \right)} \right)^2  + \left( {t_j a_j^2  + t_s a_s^2 } \right)^2 }  =  \\
 \sqrt {\left( {a_j a_s } \right)^2 \left( {2 - 2t_j t_s } \right) + a_j^4  + a_s^4  + 2t_j t_s a_s^2 a_j^2 }  = \sqrt {a_j^4  + a_s^4  + 2a_s^2 a_j^2 }  = a_j^2  + a_s^2 . \\
 \end{array}
\]
Similarly
\[
\begin{array}{l}
 \left( {a_1 ,a_{2,} ...,a_n , - it_1 a_1 , - it_2 a_2 ,..., - it_n a_n } \right) = \prod\limits_{j = 2}^n {\prod\limits_{s = 1}^{j - 1} {\left( {a_j  - a_s } \right)} } \prod\limits_{j = 2}^n {\prod\limits_{s = 1}^{j - 1} {\left( { - it_j a_j  + it_s a_s } \right)} }  \\
 \prod\limits_{j > s} {\left( { - it_j a_j  - a_s } \right)\left( { - it_s a_s  - a_j } \right)} \prod\limits_{j = 1}^n {\left( { - it_j a_j  - a_j } \right)}  =  \\
  = \left( { - 1} \right)^{\frac{{n\left( {n - 1} \right)}}{2}} i^{\frac{{n\left( {n - 1} \right)}}{2}} \left( { - 1} \right)^n \prod\limits_{j = 1}^n {a_j } \prod\limits_{j = 2}^n {\prod\limits_{s = 1}^{j - 1} {\left( {a_j  - a_s } \right)} } \prod\limits_{j = 2}^n {\prod\limits_{s = 1}^{j - 1} {\left( {t_j a_j  - t_s a_s } \right)} }  \\
 \prod\limits_{j = 1}^n {\left( {1 + it_j } \right)} \prod\limits_{j > s} {\left\{ {a_j a_s \left( {1 - t_j t_s } \right) + i\left( {t_j a_j^2  + t_s a_s^2 } \right)} \right\}}= \\
 \end{array}
\]
\[
\begin{array}{l}
  = \left( { - 1} \right)^{\frac{{n\left( {n - 1} \right)}}{2}} i^{\frac{{n\left( {n - 1} \right)}}{2}} \left( { - 1} \right)^n \sqrt 2 ^n e^{\frac{{i\pi \left( {t_1  + t_2  + ... + t_n } \right)}}{4}} \prod\limits_{j = 1}^n {a_j } \prod\limits_{j = 2}^n {\prod\limits_{s = 1}^{j - 1} {\left( {a_j  - a_s } \right)} } \prod\limits_{j = 2}^n {\prod\limits_{s = 1}^{j - 1} {\left( {t_j a_j  - t_s a_s } \right)} }  \cdot  \\
 \prod\limits_{j > s} {\left( {\left( {a_j^2  + a_s^2 } \right)e^{i\alpha _{js} } } \right)} . \\
 \end{array}
\]

Further, putting the above calculations in (18), we obtain
\[
i^{\frac{{n\left( {n + 1} \right)}}{2}} \left( { - 1} \right)^n
\sqrt 2 ^n \prod\limits_{j = 1}^n {a_j } \prod\limits_{j = 2}^n
{\prod\limits_{s = 1}^{j - 1} {\left( {a_j  - a_s } \right)} }
\prod\limits_{j = 2}^n {\prod\limits_{s = 1}^{j - 1} {\left( {t_j
a_j  - t_s a_s } \right)} } \prod\limits_{j > s} {\left( {a_j^2  +
a_s^2 } \right)}  \cdot
\]
\[
 \cdot \left( {e^{ - i\left( {\frac{{\pi \sum\limits_{j = 1}^n {t_j } }}{4} + \pi k\sum\limits_{j = 1}^n {a_j t_j }  + \sum\limits_{j > s} {\alpha _{js} } } \right)} } \right. + \left. {\left( { - 1} \right)^{\frac{{n\left( {n + 1} \right)}}{2}} e^{i\left( {\frac{{\pi \sum\limits_{j = 1}^n {t_j } }}{4} + \pi k\sum\limits_{j = 1}^n {a_j t_j }  + \sum\limits_{j > s} {\alpha _{js} } } \right)} } \right).
\]
So, if $n = 4m,\,n = 4m + 3,$ then
\[
\Delta _2  = C_2 \sum\limits_{t_j  =  \pm 1} {\left( { - 1}
\right)^l \prod\limits_{j = 2}^n {\prod\limits_{s = 1}^{j - 1}
{\left( {t_j a_j  - t_s a_s } \right)} } \cos \left( {\pi \left(
{k\sum\limits_{j = 1}^n {a_j t_j }  + \frac{{\sum\limits_{j = 1}^n
{t_j } }}{4}} \right) + \sum\limits_{j > s} {\alpha _{js} } }
\right)} ,
\]
if $n = 4m+1,\,n = 4m + 2,$  then
\[
\Delta _2  = C_3 \sum\limits_{t_j  =  \pm 1} {\left( { - 1}
\right)^l \prod\limits_{j = 2}^n {\prod\limits_{s = 1}^{j - 1}
{\left( {t_j a_j  - t_s a_s } \right)} } \sin \left( {\pi \left(
{k\sum\limits_{j = 1}^n {a_j t_j }  + \frac{{\sum\limits_{j = 1}^n
{t_j } }}{4}} \right) + \sum\limits_{j > s} {\alpha _{js} } }
\right)} ,
\]

where $C_2,\,C_3-$ are independent from $k$ and $t_ {j}$.

Further, we assume that ${a_j},\,j = 1, ...,n$ such that for all numbers $k$, with the possible exception of a finite number, exist positive constants $M,\gamma$ that the relation holds
\[\left| {{\Delta _2}} \right| > \frac{M}{{{k^\gamma }}} > 0,\eqno(19)\]
the fulfillment of inequality (19) was shown above, at least for ${a_j} \in N$ and ${a_j} \in Q$, provided that (17) is satisfied, with the exponent $\gamma = 0.$ Next, not to increase the number of notations; all positive constants will be denoted by $ M. $

We turn to the condition when there is a solution to the Dirichlet problem. We will seek a solution in the form
\[u\left( {x,y} \right) = \sqrt 2 \sum\limits_{k = 1}^\infty  {{u_k}\left( y \right)\sin \pi kx} ,\eqno(20)\]
the theorem holds.

\textbf{Theorem 3.} Let coefficients ${a_j},\,j = 1,...,n$ such that condition (19) is fulfilled for them and, moreover, the relations hold
\[{\varphi _j}\left( x \right),{\psi _j}\left( x \right) \in {C^{\left( {2n + 1 + \alpha } \right)}}\left[ {0,1} \right],\,\,j = 1,...,n,\]
\[\left. {\left( {\varphi _j^{\left( {s - 1} \right)}\left( x \right)\sin \left( {\pi kx + \frac{\pi }{2}s} \right)} \right)} \right|_{x = 0}^{x = 1} = 0,\,\left. {\left( {\psi _j^{\left( {s - 1} \right)}\left( x \right)\sin \left( {\pi kx + \frac{\pi }{2}s} \right)} \right)} \right|_{x = 0}^{x = 1} = 0,\]
\[j = 1,...,n,\,s = 1,...,2n + 1 + \alpha ,\]
\[\alpha  = \left\{ \begin{array}{l}
\gamma ,\,\,\,\gamma  \in N \cup \left\{ 0 \right\},\\
\left[ \gamma  \right] + 1,\,\,0 < \gamma  \notin N.
\end{array} \right.\]
Then series (20) is a classical solution to the Dirichlet problem.

\textbf{Proof.} Let $y \ne 0.$ Formally, we differentiate expression (20) $2n$ times in the variable
\[\frac{{{\partial ^{2n}}u\left( {x,y} \right)}}{{\partial {x^{2n}}}} = \sqrt 2 \sum\limits_{k = 1}^\infty  {{{\left( { - 1} \right)}^n}{{\left( {\pi k} \right)}^{2n}}{u_k}\left( y \right)\sin \pi kx} ,\]
so we have
\[\left| {\frac{{{\partial ^{2n}}u\left( {x,y} \right)}}{{\partial {x^{2n}}}}} \right| \le M\sum\limits_{k = 1}^\infty  {{k^{2n}}\left| {{u_k}\left( y \right)} \right|}  \le M\sum\limits_{k = 1}^\infty  {{k^{2n}}\sum\limits_{s = 0}^{n - 1} {\frac{{\left( {\left| {{\varphi _{sk}}} \right| + \left| {{\psi _{sk}}} \right|} \right)}}{{\left| {{\Delta _2} + {\Delta _k}} \right|}}} } ,\]

we show the convergence of the series at $s = 0$; for the remaining terms, the convergence is shown identically. Note that the following cases are possible:

 1. As $\Delta _k$ tends to zero exponentially, then for sufficiently large values of $a_i$,$ i = 1, ..., n$, the relation  ${\left| {{\Delta _2} + {\Delta _k}} \right|}>0$.  In this case we will have

\[M\sum\limits_{k = 1}^\infty  {{k^{2n}}\frac{{\left| {{\varphi _{0k}}} \right|}}{{\left| {{\Delta _2} + {\Delta _k}} \right|}}}  \le M\sqrt {\sum\limits_{k = 1}^\infty  {\frac{1}{{{k^2}}}} } \sqrt {\sum\limits_{k = 1}^\infty  {{{\left( {\frac{{\left| {{\varphi _{0k}}} \right|{k^{2n + 1}}}}{{\left| {{\Delta _2} + {\Delta _k}} \right|}}} \right)}^2}} }  \le M\frac{\pi }{{\sqrt 6 }}\sqrt {\sum\limits_{k = 1}^\infty  {{{\left( {\left| {{\varphi _{0k}}} \right|{k^{2n + 1 + \alpha }}} \right)}^2}} }  = \]
\[ = M\frac{\pi }{{\sqrt 6 }}{\left\| {\varphi _0^{\left( {2n + 1 + \alpha } \right)}\left( x \right)} \right\|_{{L_2}\left[ {0,1} \right]}}\]
(Parseval equality was used here). Similarly, the possibility of termwise differentiation $2n$ times with respect to the variable $y$ is shown.

\textbf{Conclusion}: Under the conditions of Theorem 3 and sufficiently large values of $a_i$ - there is a unique classical solution to the Cauchy problem.

2.It may turn out that, even if condition (19) is satisfied, the expression $\Delta = 0,$ for some finite values ??$k = {k_1}, {k_2}, ..., {k_p} <{k_0}, $ where $ {k_1} <{k_2} <... <{k_p} $, and $p$ are given natural numbers. Then, for the solvability of the Dirichlet problem, it suffices to satisfy the conditions ${\varphi _{sk}} = \sqrt 2 \int\limits_0^1 {{\varphi _s}\left( x \right)\sin \pi kdx}  = 0,$ ${\psi _{sk}} = \sqrt 2 \int\limits_0^1 {{\psi _s}\left( x \right)\sin \pi kxdx} ,\,s = \overline {0,(n - 1)} ,k = {k_1},{k_2},...,{k_p},$ and the solution itself will have the form
\[u\left( {x,y} \right) = \sqrt 2 \sum\limits_{\mathop {k = 1}\limits_{k \ne {k_1},{k_2},...,{k_p}} }^\infty  {{u_k}\left( y \right)\sin \pi kx}  + \sqrt 2 \sum\limits_{i = 1}^p {\widetilde {{u_{{k_i}}}}\left( y \right)\sin \pi {k_i}x} ,\]
Here function $\widetilde {{u_{{k_i}}}}\left( y \right)$ - is a nonzero solution to a homogeneous system(9).

\textbf{Conclusion:}
In general, when the conditions of Theorem 3 are satisfied, a solution exists, but uniqueness may be violated.

\textbf{Theorem is proved. }\\

\begin{center}
\textbf{Reference}
  \end{center}
1. Sneddon I.N., Berry D.S. (1958) The Classical Theory of Elasticity. In: Flugge S. (eds) Elasticity and Plasticity / Elastizitat und Plastizitat. Handbuch der Physik / Encyclopedia of Physics, vol 3 / 6. Springer, Berlin, Heidelberg \\
2. Umezawa H. Quantum Field Theory.North-Holland Pub., 1956. — 364 p.\\
3. Biсadze A.V. Selected Works.-Nalchik: Publishing House Science.-2012.p.400 [In Russian].\\
4. Smirnov, M. M. The first boundary value problem for a certain equation of mixed type. (Russian) Differencialnye Uravnenija 4 1968 1663–1666.\\
5. Smirnov, M. M. A boundary value problem with shift for a fourth order equation of mixed-composite type. (Russian) Differencialnye Uravnenija 11 (1975), no. 9, 1678–1686, 1718.\\
6. Zegalov, V. I. Boundary value problem for a mixed-type equation of higher order. Dokl. Akad. Nauk SSSR 136 274–276 (Russian); translated as Soviet Math. Dokl. 2 1961 45–47.\\
7. Sabitov K. B.On positiveness of a solution to inhomogeneous mixed type equation of higher order. Russian Math. (Iz. VUZ), 60:3 (2016), 56–62\\
8. Sabitov, K. B. Solutions of fixed sign of a higher-order inhomogeneous equation of mixed elliptic-hyperbolic type. (Russian) Mat. Zametki 100 (2016), no. 3, 433–440; translation in Math. Notes 100 (2016), no. 3-4, 458–464\\
9. Zegalov, V. I. The equation with constant coefficients that is linear with respect to the Lavrent?ev-Bicadze operator. (Russian) Trudy Sem. Kraev. Zadacam Vyp. 6 (1969), 44–51.\\
10.Ptashnik B.I. Incorrect boundary value problems for differential
partial differential equations // Kiev, Naukova Dumka, 1984, -
p.264 [In Russian]. \\
11. Sabitov, K. B. The Dirichlet problem for higher-order partial differential equations. (Russian) Mat. Zametki 97 (2015), no. 2, 262–276; translation in Math. Notes 97 (2015), no. 1-2, 255–267\\
12. Irgashev B. Yu. On one boundary-value problem for an equation of higher even order. Russian Math. (Iz. VUZ), 61:9 (2017), 10–26\\
13. Irgashev  B. Yu. On spectral problem for one equation of high even order. Russian Math. (Iz. VUZ), 60:7 (2016), 37–46\\
14. Bicadze, A. V. Incorrectness of Dirichlet's problem for the mixed type of equations in mixed regions. (Russian) Dokl. Akad. Nauk SSSR 122 1958 167–170.\\
15. Sabitov, K. B.; Khadzhi, I. A. A boundary value problem for the Lavrent?ev-Bitsadze equation with an unknown right-hand side. (Russian) Izv. Vyssh. Uchebn. Zaved. Mat. 2011, no. 5, 44–52; translation in Russian Math. (Iz. VUZ) 55 (2011), no. 5, 35–42\\
16.Khadzhi, I.A. Inverse problem for equations of mixed type with Lavrent’ev-Bitsadze operator. Math Notes 91, 857–867 (2012). https://doi.org/10.1134/S0001434612050331

\end{document}